\documentclass[oneside,english]{amsart}
\usepackage[T1]{fontenc}
\usepackage[latin9]{inputenc}
\usepackage{amsthm}
\usepackage{amstext}
\usepackage{amssymb}
\usepackage{amsxtra}
\usepackage{esint}
\usepackage{graphicx}
\usepackage{latexsym}

\makeatletter
\numberwithin{equation}{section}
\numberwithin{figure}{section}
\theoremstyle{plain}
\newtheorem{thm}{\protect\theoremname}[section]
\newtheorem{prop}[thm]{\protect\propositionname}
\theoremstyle{plain}
\theoremstyle{definition}

\theoremstyle{plain}

\theoremstyle{plain}
\newtheorem{rem}[thm]{\protect\remarkname}
\theoremstyle{plain}
\makeatother

\usepackage{babel}
\usepackage{color}
\providecommand{\definitionname}{Definition}
\providecommand{\lemmaname}{Lemma}
\providecommand{\theoremname}{Theorem}
\providecommand{\corollaryname}{Corollary}
\providecommand{\remarkname}{Remark}
\providecommand{\propositionname}{Proposition}

\DeclareMathOperator{\loc}{loc}

\DeclareMathOperator{\cp}{cap}

\begin{document}

\title[On the theory of generalized quasiconformal mappings]
{On the theory of generalized quasiconformal mappings}

\author{Vladimir Gol'dshtein, Evgeny Sevost'yanov, Alexander Ukhlov}
\begin{abstract}
We study generalized quasiconformal mappings in the context of the inverse Poletsky inequality. We consider the local behavior and the  boundary behavior of mappings with the inverse Poletsky inequality. In particular, we obtain logarithmic H\"{o}lder
continuity for such classes of mappings.
\end{abstract}
\maketitle

\section{Introduction }

This article is devoted to the H\"older continuity of generalized quasiconformal mappings $f:D\to D'$ which are defined by capacity
(moduli) inequalities. The method of capacity (moduli) inequalities arises
to the Gr\"{o}tzsch problem and was introduced in \cite{A66}. In subsequent works
(see, for example, \cite{MRV, MRSY} and \cite{Ri}) the conformal modulus
method was used in the theory of quasiconformal (quasiregular) mappings and its generalizations.
The  classes of mappings generating bounded composition operators on Sobolev
spaces  \cite{VU04,VU05} arise in the geometric analysis of PDE \cite{GS82,M69}.
These mappings are called weak $(p,q)$-quasiconformal mappings \cite{GGR95,VU98}
and can be characterized by the inverse capacitory (moduli) Poletsky inequality \cite{U93}
$$
{\rm cap}_q^{1/q}(f^{-1}(E),f^{-1}(F);D) \leqslant
K_{p,q}(\varphi;\Omega) {\rm cap}_{p}^{1/p}(E,F;D'),\,\,1<q\leqslant
p<\infty.
$$
The detailed study of the mappings with the inverse conformal Poletsky
inequality for  modulus of paths was given in \cite{SSD}, \cite{SevSkv} and~\cite{Sev$_3$}.
In this case $p=q=n$ the H\"{o}lder continuity, the continuous boundary extension, and the
behavior on the closure of domains were obtained.

In the recent works \cite{GSU23,MU23} were considered connections between weak
$(p,q)$-quasiconformal mappings and $Q$-homeomorphisms.
In the present article we suggest an approach to the generalized quasiconformal mappings
which is based on the following integral inequality
\begin{equation*}
\int\limits _{D}|\nabla(u\circ f(x))|^{q}~dm(x)\leqslant\int\limits
_{D'}|\nabla u(y)|^{q}Q_q(y)~dm(y),\,\,u\in C^1(D').
\end{equation*}
Depending on the properties of the function $Q_q$ we obtain various classes
of the generalized quasiconformal mappings:
$BMO$-quasiconformal mappings, weak $(p,q)$-quasiconformal mappings, $Q$-mappings and so on.

The weak $(p,q)$-quasiconformal mappings have significant
applications in the spectral theory of elliptic operators \cite{GU16,GU17}. The H\"{o}lder continuity of weak $(p,q)$-quasiconformal mappings was considered in \cite{VU98}.
In the recent article \cite{GSU} was considered the boundary behavior of the weak
$(p,q)$-quasiconformal mappings. In this article we study the logarithmic H\"{o}lder continuity,
the continuous boundary extension, and the behavior in the closure of domains of
non-homeomorphic generalizations of quasiconformal mappings.

\vskip 0.2cm
Let us give the basic definitions. Let $\Gamma$ be a family of paths
$\gamma$ in ${\mathbb R}^n$. A Borel function $\rho:{\mathbb R}^n\,\rightarrow [0,\infty] $ is called {\it admissible} for
$\Gamma$ if
\begin{equation}\label{eq1.4}
\int\limits_{\gamma}\rho (x) |dx|\geqslant 1
\end{equation}
for all (locally rectifiable) paths $\gamma \in \Gamma$. In this
case, we write: $\rho \in {\rm adm} \,\Gamma$. Given a number
$q\geqslant 1,$ {\it $q$-modulus} of the family of paths $\Gamma$ is
defined as
\begin{equation}\label{eq1.3gl0} M_q(\Gamma)=\inf\limits_{\rho \in \,{\rm adm}\,\Gamma}
\int\limits_D \rho^{\,q} (x)\, dm(x)\,.
\end{equation}

Let $x_0\in\overline{D}$, $x_0\ne\infty$, then
\begin{equation}\label{eq1E}
B(x_0, r)=\{x\in {\mathbb R}^n: |x-x_0|<r\}\,, \quad {\mathbb B}^n=B(0,
1)\,, \end{equation}
$$S(x_0,r) = \{
x\,\in\,{\mathbb R}^n : |x-x_0|=r\}\,, S_i=S(x_0, r_i)\,,\quad
i=1,2\,,$$
\begin{equation*}\label{eq1**}
A=A(x_0, r_1, r_2)=\{ x\,\in\,{\mathbb R}^n : r_1<|x-x_0|<r_2\}\,.
\end{equation*}
Given sets $E$, $F\subset\overline{{\mathbb R}^n}$ and a domain $D\subset {\mathbb R}^n$, we denote $\Gamma(E,F,D)$ a family of all
paths $\gamma:[a,b]\rightarrow \overline{{\mathbb R}^n}$ such that $\gamma(a)\in E,\gamma(b)\in\,F $ and $\gamma(t)\in D$ for all $t \in
(a, b)$.

Let $Q:{\mathbb R}^n\rightarrow [0, \infty]$ be a Lebesgue
measurable function. We say that {\it $f$ satisfies the Poletsky
inverse inequality with respect to $q$-modulus} at a point $y_0\in
f(D),$ $1<q<\infty$, if the moduli inequality
\begin{equation}\label{eq2*A}
M_q(\Gamma(E, F, D))\leqslant
\int\limits_{A(y_0,r_1,r_2)\cap f(D)} Q(y)\cdot \eta^{\,q}(|y-y_0|)\, dm(y)
\end{equation}
holds for any continua $E\subset f^{\,-1}(\overline{B(y_0, r_1)})$,
$F\subset f^{\,-1}(f(D)\setminus B(y_0, r_2))$,
$0<r_1<r_2<r_0=\sup\limits_{y\in f(D)}|y-y_0|$, and any Lebesgue
measurable function $\eta: (r_1,r_2)\rightarrow [0,\infty ]$ such
that
\begin{equation}\label{eqA2}
\int\limits_{r_1}^{r_2}\eta(r)\, dr\geqslant 1\,.
\end{equation}

\medskip
The case $q=n$ was studied in details in~\cite{SSD}, cf.~\cite{SevSkv} and \cite{Sev$_3$}. The present article is
dedicated to the case $q\ne n$.

\medskip
Let us formulate the main results of this manuscript. Recall that a mapping $f:D\rightarrow {\mathbb R}^n$ is called {\it discrete} if a
pre-image $\{f^{-1}\left(y\right)\}$ of each point $y\,\in\,{\mathbb R}^n$ consists of isolated points, and {\it open} if the image of
any open set $U\subset D$ is an open set in ${\mathbb R}^n$.  The mapping $f$ of the domain $D$ onto $D^{\,\prime}$ is called {\it closed} if
$f(E)$ is closed in $D^{\,\prime}$ for any of the closed $E\subset D$ (see, e.g., \cite[Section~3]{Vu}).

In the extended Euclidean $n$-dimensional space $\overline{{{\mathbb R}}^n}={{\mathbb
R}}^n\cup\{\infty\}$, a {\it spherical (chordal) metric} is defined as $h(x,y)=|\pi(x)-\pi(y)|$, where $\pi$ is a stereographic
projection of $\overline{{{\mathbb R}}^n}$ onto the sphere $S^n(\frac{1}{2}e_{n+1},\frac{1}{2})$ in ${{\mathbb R}}^{n+1}$. Namely:
\begin{equation}\label{eq3C}
 h(x,y)=\frac{|x-y|}{\sqrt{1+{|x|}^2} \sqrt{1+{|y|}^2}}\,, \, \, x\ne \infty\ne y,\,\,h(x,\infty)=\frac{1}{\sqrt{1+{|x|}^2}}.
\end{equation}
(see, e.g., \cite[definition~12.1]{Va}).
Given sets $A, B\subset \overline{{\mathbb R}^n},$ we put
$$h(A, B)=\inf\limits_{x\in A, y\in B}h(x, y)\,,\quad h(A)=\sup\limits_{x, y\in A}h(x ,y)\,,$$
where $h$ is defined in~(\ref{eq3C}).
In addition, we put
$${\rm dist}\,(A, B)=
\inf\limits_{x\in A, y\in B}|x-y|\,,\quad {\rm
diam}(A)=\sup\limits_{x, y\in A}|x-y|\,.$$

Let $D\subset {\mathbb R}^n$, $n\geqslant 2$, be a domain. For a number $1\leqslant q<\infty$ and a Lebesgue measurable function $Q:{\mathbb
R}^n\rightarrow [0, \infty],$ we denote by $\mathfrak{F}^{q}_Q(D)$ a family of all open discrete mappings $f:D\rightarrow {\mathbb R}^n$
such that relation~(\ref{eq2*A}) holds for any $y_0\in f(D)$, for any continua
$$E\subset f^{\,-1}(\overline{B(y_0, r_1)}), \,\,F\subset f^{\,-1}(f(D)\setminus B(y_0, r_2)),\,\,
0<r_1<r_2<r_0=\sup\limits_{y\in f(D)}|y-y_0|,
$$
and any Lebesgue measurable function $\eta: (r_1,r_2)\rightarrow [0,\infty ]$ with condition~(\ref{eqA2}).

The following theorem holds.

\begin{thm} \label{th1}
{\sl Let $f\in \mathfrak{F}^{q}_Q({\mathbb B}^n)$, $q\geqslant n$. Suppose that $Q\in L^1({\mathbb
R}^n)$ and $K$ is a compact set in ${\mathbb B}^n$. Then the inequality
\begin{equation}\label{eq2C}
|f(x)-f(y)|\leqslant C_n \cdot\frac{(\Vert Q\Vert_1)^{\frac{1}{q}}}{\log^{\frac{1}{n}}\left(1+\frac{r_0}{2|x-y|}\right)}, \,\,r_0 =d(K, \partial {\mathbb B}^n),
\end{equation}
holds for all $x, y\in K$, where $\Vert Q\Vert_1$ denotes the $L^1$-norm of the function $Q$ in
${\mathbb R}^n$ and a constant $C_n>0$ depends on $n$ and $q$ only.}
\end{thm}

Let $D\subset {\mathbb R}^n$ be a domain. Then $D$ is called {\it locally connected at the point} $x_0\in\partial D$, if for any neighborhood $U$ of $x_0$ there is a neighborhood $V\subset U$ of this point such that $V\cap D$ is connected. The domain $D$ is locally connected on $\partial D$, if $D$ is locally connected at every point $x_0\in\partial D$.
The domain $D\subset {\mathbb R}^n$ is called {\it finitely connected at
the point} $x_0\in\partial D,$ if for any neighborhood $U$ of $x_0$
there is a neighborhood $V\subset U$ of this point such that the set
$V\cap D$ consists of a finite number of components (see, e.g.,
\cite{Vu}). The domain $D$ is finitely connected on $\partial D,$ if
$D$ is finitely connected at every point $x_0\in\partial D.$

Let $\partial D$ be a boundary of the domain $D\subset\mathbb R^n$. Then the boundary $\partial D$ is called {\it weakly flat} at the
point $x_0\in \partial D,$ if for each $P>0$ and for any neighborhood $U$ of this point there is a neighborhood $V\subset U$
of the same point such that $M(\Gamma(E, F, D))>P$ for any continua $E, F\subset D$ that intersect $\partial U$ and $\partial V$. The
boundary of a domain $D$ is called weakly flat if the corresponding property holds at any point of $\partial D$.

\medskip
Let $D, D^{\,\prime}$ be domains in $\mathbb R^n$. For given numbers $n\leqslant q<\infty$, $\delta>0$, a continuum $A\subset
D^{\,\prime}$ and an arbitrary Lebesgue measurable function $Q:D^{\,\prime}\rightarrow [0, \infty],$ we denote by ${\mathfrak
S}^{q}_{\delta, A, Q }(D, D^{\,\prime})$ a family of all open discrete and closed mappings $f$ of $D$ onto $D^{\,\prime}$
satisfying the condition~(\ref{eq2*A}) for any $y_0\in D^{\,\prime},$ any compacts
$$
E\subset f^{\,-1}(\overline{B(y_0,r_1)}),\,\,F\subset f^{\,-1}(D^{\,\prime}\setminus B(y_0, r_2)),\,\,0<r_1<r_2<r_0=\sup\limits_{y\in D^{\,\prime}}|y-y_0|,
$$
and any Lebesgue measurable function $\eta: (r_1,r_2)\rightarrow [0,\infty
]$ with the condition~(\ref{eqA2}), such that $h(f^{\,-1}(A),
\partial D)\geqslant~\delta.$ The following statement holds.

\medskip
\begin{thm}\label{th2}
{\sl\,Let $D\subset\mathbb R^n$ be a bounded with a weakly flat boundary. Suppose that, for any point $y_0\in
\overline{D^{\,\prime}}$ and $0<r_1<r_2<r_0:=\sup\limits_{y\in
D^{\,\prime}}|y-y_0|$ there is a set $E\subset[r_1, r_2]$ of a
positive linear Lebesgue measure such that the function $Q$ is
integrable on $S(y_0, r)$ for every $r\in E.$ If $D^{\,\prime}$ is
locally connected on its boundary, then any $f\in{\mathfrak
S}^{q}_{\delta, A, Q }(D, D^{\,\prime})$ has a continuous extension
$\overline{f}:\overline{D}\rightarrow \overline{D^{\,\prime}},$
$\overline{f}(\overline{D})=\overline{D^{\,\prime}},$ and the family
${\mathfrak S}^{q}_{\delta, A, Q }(\overline{D},
\overline{D^{\,\prime}}),$ which consists of all extended mappings
$\overline{f}:\overline{D}\rightarrow \overline{D^{\,\prime}},$ is
equicontinuous in $\overline{D}.$

In particular, the statement of Theorem~\ref{th2} is fulfilled if
the above condition on $Q$ is replaced by a simpler one: $Q\in
L^1(D^{\,\prime}).$
 }
\end{thm}

\medskip
\begin{rem}
\label{rem1} In Theorem~\ref{th2}, the equicontinuity must be
understood with respect to the Euclidean metric in the preimage
under the mapping, and the chordal metric in the image, i.e., for
any $\varepsilon>0$ there is $\delta=\delta(\varepsilon, x_0)>0 $
such that the condition $|x-x_0|<\delta,$ $x\in D,$ implies that
inequality $h(\overline{f}(x, \overline{f}(x_0))<\varepsilon$ holds
for any $\overline{f}\in{\mathfrak S}^{q}_{\delta, A, Q }(\overline{D},
\overline{D^{\,\prime}}).$
\end{rem}

\section{On the integral inverse Poletsky inequality}

In this section we suggest an approach to the generalized quasiconformal mappings which is based on the following integral inequality
\begin{equation*}
\int\limits _{D}|\nabla(u\circ f(x))|^{q}~dm(x)\leqslant\int\limits
_{D'}|\nabla u(y)|^{q}Q_q(y)~dm(y),\,\,u\in C^1(D').
\end{equation*}
This approach allows to unify various generalizations of quasiconformal mappings, such as mappings which generate bounded composition operators on  seminormed Sobolev spaces and $Q$-mappings. We explain that both concepts of generalizations are very close one to another and, in some sense, represent similar classes. Of course, it is a subject of more deep study. We are trying to put attention of readers to this useful interplay.

Let $D$ be a domain in the Euclidean space $\mathbb R^n$, $n\geqslant 2$. The Sobolev space
$W^1_p(D)$, $1\leqslant p\leqslant\infty$,
is defined as a Banach space of locally integrable weakly differentiable functions
$u:D\to\mathbb{R}$ equipped with the following norm:
\[
\|u\mid W^1_p(D)\|=\| u\mid L_p(D)\|+\|\nabla u\mid L_p(D)\|,
\]
where $\nabla u$ is the weak gradient of the function $u$.

The homogeneous seminormed Sobolev space $L^1_p(D)$, $1\leqslant p\leqslant\infty$, is defined as a space
of locally integrable weakly differentiable functions $u:D\to\mathbb{R}$ equipped
with the following seminorm:
\[
\|u\mid L^1_p(D\|=\|\nabla u\mid L_p(D)\|.
\]

In accordance with the non-linear potential theory \cite{MH72}
we consider elements of Sobolev spaces $W^1_p(\Omega)$ as equivalence
classes up to a set of $p$-capacity zero  \cite{M}.

Suppose $f:D\to\mathbb{R}^{n}$ is a mapping of the Sobolev class $W^1_{1,\loc}(D;\mathbb R^n)$.
Then the formal Jacobi
matrix $D f(x)$ and its determinant (Jacobian) $J(x,f)$
are well defined at almost all points $x\in D$. The norm $|D f(x)|$ is the operator norm of $D f(x)$.

Recall the change of variable formula for the Lebesgue integral \cite{H93}.
Let a mapping $f : D\to \mathbb R^n$ belongs to $W^1_{1,\loc}(D;\mathbb R^n)$.
Then there exists a measurable set $S\subset D$, $|S|=0$ such that  the mapping $f:D\setminus S \to \mathbb R^n$ has the Luzin $N$-property and the change of variable formula
\begin{equation}
\label{chvf} \int\limits_E u\circ f (x)
|J(x,f)|~dm(x)=\int\limits_{\mathbb R^n\setminus \varphi(S)}
u(y)N_f(E,y)~dm(y)
\end{equation}
holds for every measurable set $E\subset D$ and every non-negative measurable function $u: \mathbb R^n\to\mathbb R$. Here

Now let $D$ and $D'$ be domains in Euclidean space $\mathbb{R}^{n}$,
$n\geqslant2$. We consider a homeomorphism $f:D\to D'$ of the class
$W^1_{1,\loc}(D;D')$ which has finite distortion. Recall that the
mapping $ f$ is called the mapping of finite distortion if
$|Df(x)|=0$ for almost all $x\in Z=\{z\in D:J(x,f)=0\}$.

By using the composition of functions $u\in C^{1}(D)$ with this homeomorphism $f:D\to D'$ we obtain the following inequality
\begin{multline*}
\|u\circ f\mid L_{q}^{1}(D)\|^q:=\int\limits _{D}|\nabla(u\circ f(x))|^{q}~
dm(x)\leqslant\int\limits _{D}|\nabla u(f(x))|^{q}|Df(x))|^{q}~dm(x)\\
=\int\limits _{D\setminus Z}|\nabla
u(f(x))|^{q}|J(x,f)||Df(x))|^{q}|J(x,f)|^{-1}~dm(x).
\end{multline*}
By the change of variables formula \cite{H93} we have the following {\it integral inverse Poletsky inequality}
\begin{equation}
\label{baseeq} \int\limits _{D}|\nabla(u\circ f(x))|^{q}~dm(x)\leqslant\int\limits _{D'}|\nabla
u(y)|^{q}Q_q(y)~dm(y),
\end{equation}
where
\[
Q_{q}(y):=\begin{cases}
\frac{|Df(x)|^{q}}{|J(x,f)|},\,\, & x=f^{-1}(y)\in\Omega\setminus(S\cup Z),\\
\,\,0,\,\, & x=f^{-1}(y)\in S\cup Z.
\end{cases}
\]
The characterization of mappings which generate bounded composition operators on Sobolev spaces in terms of integrability of this distortion function $Q_q$ was given in \cite{VU02} (see, also, \cite{VU04,VU05}).

Depending on the properties of the distortion function $Q_{q}(y)$
we obtain different classes of generalized quasiconformal
mappings. Let us recall the notion of the variational $p$-capacity \cite{GResh}.
The condenser in the domain $D\subset \mathbb R^n$ is the pair $(E,F)$
of connected closed relatively to $D$ sets $E,F\subset D$.
Recall that a continuous function $u\in L_p^1(D)$ is called an
admissible function for the condenser $(E,F)$, denoted $u\in W_0(E, F)$, if the set
$E\cap D$ is contained in some connected component of the set ${\rm Int}\{x: u(x)=0\}$, the set $F\cap D$ is contained in some to the connected
component of the set ${\rm Int}\{x: u(x)=1\}$. Then we call as a $p$-capacity of the condenser
$(E,F)$ relatively to a domain $D$
the value
$$
{{\cp}}_p(E,F;\Omega)=\inf\|u\vert L_p^1(D)\|^p,
$$
where the greatest lower bond is taken over all admissible for the condenser $(E,F)\subset D$
functions. If the condenser have no admissible functions we put the capacity is equal to infinity.

\vskip 0.2cm
\noindent
{\it The case of $K$-quasiconformal mappings.} Let $q=n$ and
    $$
    {\rm ess}\sup_{y\in D'} Q_n(y)={\rm ess}\sup_{y\in D'}\frac{|Df(f^{-1}(y))|^n}{|J(f^{-1}(y),f)|}=K_n<\infty.
    $$
Then by the inequality (\ref{baseeq}) for any condenser $(E,F)\subset D'$ the inequality
$$
\cp_n(f^{-1}(E),f^{-1}(F);D)\leqslant K_n \cp_n(E,F;D')
$$
holds. Hence $f$ is a $K_n$-quasiconformal mapping \cite{Va}. From another side quasiconformal mappings generate bounded composition operators on Sobolev spaces $L^1_n(D')$ and $L^1_n(D)$ \cite{VG75}.

The special case represent conformal mappings that corresponds to the case $q=n=2$ and $K=1$. In this case we have isometries of Sobolev spaces $L^1_2(D')$ and $L^1_2(D)$.

\vskip 0.2cm
\noindent
{\it The case of $q$-quasiconformal mappings.} Let $1<q<\infty$ and
    $$
    {\rm ess}\sup_{y\in D'} Q_q(y)={\rm ess}
    \sup_{y\in D'}\frac{|Df(f^{-1}(y))|^q}{|J(f^{-1}(y),f)|}=K_q<\infty.
    $$
Then by the inequality (\ref{baseeq}) for any condenser $(E,F)\subset D'$ the inequality
$$
\cp_q(f^{-1}(E),f^{-1}(F);D)\leqslant K_q \cp_q(E,F;D')
$$
holds. Hence $f$ is a $q$-quasiconformal mapping \cite{VU98}. From another side by \cite{GGR95,VU98} $q$-quasiconformal mappings generate bounded composition operators on Sobolev spaces $L^1_q(D')$ and $L^1_q(D)$.

\vskip 0.2cm
\noindent
{\it The case of $(p,q)$-quasiconformal mappings.} Let $1<q<p<\infty$
and $Q_q\in L_s(\Omega)$, $s>1$. Then by the H\"older inequality
\begin{multline*}
\left(\int\limits _{D}|\nabla(u\circ
f)|^{q}~dm(x)\right)^{\frac{1}{q}}
\leqslant\left(\int\limits _{D'}|\nabla u(y)|^{q}Q_q(y)~dm(y)\right)^{\frac{1}{q}}\\
\leqslant \left(\int\limits_{D'}Q_q^s(y)~dm(y)\right)^{\frac{1}{qs}}
\left(\int\limits_{D'} |\nabla u(y)|^{q\frac{s}{s-1}}
~dm(y)\right)^{\frac{s-1}{qs}}.
\end{multline*}
Denote $p=q\frac{s}{s-1}$. Then $s=p/(p-q)$ and we obtain
$$
\left(\int\limits _{D}|\nabla(u\circ
f)|^{q}~dm(x)\right)^{\frac{1}{q}}\leqslant
\left(\int\limits_{D'}Q_q^{\frac{p}{p-q}}(y)~dm(y)\right)^{\frac{p-q}{pq}}
\left(\int\limits_{D'}|\nabla u(y)|^{p} ~dm(y)\right)^{\frac{1}{p}}.
$$
Hence \cite{U93} for any condenser $(E,F)\subset D'$ the inequality
$$
\cp^{\frac{1}{q}}_q(f^{-1}(E),f^{-1}(F);D)\leqslant
\left(\Phi(D'\setminus(E\cup F))\right)^{\frac{p-q}{pq}}
\cp^{\frac{1}{p}}_p(E,F;D')
$$
holds, where
$$
\Phi(D'\setminus(E\cup F))=\int\limits_{D'\setminus(E\cup
F)}Q_q^{\frac{p}{p-q}}(y)~dm(y).
$$
So $f$ is a $(p,q)$-quasiconformal mapping \cite{U93}.
From another side by \cite{U93} $(p,q)$-quasiconformal
mappings generate bounded composition operators on Sobolev spaces $L^1_p(D')$ and $L^1_q(D)$.

\vskip 0.2cm
\noindent
{\it The case of $Q$-mappings.} Let $q=n$ and $Q_q\in L_1(D)$.
Then we have the class of mappings with capacitory inverse Poletsky
inequality which was intensively studied recently \cite{SSD}, \cite{SevSkv} and~\cite{Sev$_3$}.

\vskip 0.2cm
So we can conclude that the integral inequality
$$
\int\limits_{D}|\nabla (u\circ\ f)|^q~dm(x)\leqslant
\int\limits_{D'}|\nabla u(y)|^q Q (y)~dm(y)
$$
is the basic tool for generalizations of quasiconformal mappings. In the present
work we consider the H\"{o}lder continuity and the continuous boundary extension of continuous mappings $f: D:\mathbb R^n$ in the case $q\ne n$ and $Q_q\in L_1(D)$.
This class of mappings derives properties of mappings $(p,q)$-quasiconformal
mappings which are important in the spectral theory of elliptic operators.

In the case of connected closed relatively to $D$ sets $E,F\subset D$ the notions of the capacity and the modulus coincide, but in view of suggested techniques we will use the notion of the modulus.

\section{On the H\"{o}lder continuity of mappings}

Let us first formulate the important topological statement, which is repeatedly used later (see, for example,
\cite[theorem~1.I.5.46]{Ku}).

\medskip
\begin{prop}\label{pr2}
{\sl\, Let $A$ be a set in a topological space $X.$ If the set $C$ is connected,  $C\cap A\ne \varnothing$ and $C\setminus A \ne \varnothing$, then
$C\cap
\partial A\ne\varnothing$.}
\end{prop}

\medskip
Let $D\subset {\mathbb R}^n,$ $f:D\rightarrow {\mathbb R}^n$ be a discrete
open mapping, $\beta: [a,\,b)\rightarrow {\mathbb R}^n$ be a path, and
$x\in\,f^{\,-1}(\beta(a)).$ A path $\alpha: [a,\,c)\rightarrow D$ is
called a {\it maximal $f$-lifting} of $\beta$ starting at $x,$ if
$(1)\quad \alpha(a)=x\,;$ $(2)\quad f\circ\alpha=\beta|_{[a,\,c)};$
$(3)$\quad for $c<c^{\prime}\leqslant b,$ there is no a path
$\alpha^{\prime}: [a,\,c^{\prime})\rightarrow D$ such that
$\alpha=\alpha^{\prime}|_{[a,\,c)}$ and $f\circ
\alpha^{\,\prime}=\beta|_{[a,\,c^{\prime})}.$ Similarly, we may
define a maximal $f$-lifting $\alpha: (c,\,b]\rightarrow D$ of a
path $\beta: (a,\,b]\rightarrow {\mathbb R}^n$ ending at
$x\in\,f^{\,-1}(\beta(b)).$ The maximal lifting $\alpha:
[a,\,c)\rightarrow D$ of the path $\beta: [a,\,b)\rightarrow {\mathbb
R}^n$ at the mapping $f$ with the origin at the point $x $ is called
{\it whole (total)} if, in the above definition, $c=b.$ The
following assertion holds (see~\cite[Lemma~3.12]{MRV$_3$}).

\medskip
\begin{prop}\label{pr3}
{\sl Let $f:D\rightarrow {\mathbb R}^n,$ $n\geqslant 2,$ be an open
discrete mapping, let $x_0\in D,$ and let $\beta: [a,\,b)\rightarrow
{\mathbb R}^n$ be a path such that $\beta(a)=f(x_0)$ and such that
either $\lim\limits_{t\rightarrow b}\beta(t)$ exists, or
$\beta(t)\rightarrow \partial f(D)$ as $t\rightarrow b.$ Then
$\beta$ has a maximal $f$-lifting $\alpha: [a,\,c)\rightarrow D$
starting at $x_0.$ If $\alpha(t)\rightarrow x_1\in D$ as
$t\rightarrow c,$ then $c=b$ and $f(x_1)=\lim\limits_{t\rightarrow
b}\beta(t).$ Otherwise $\alpha(t)\rightarrow \partial D$ as
$t\rightarrow c.$}
\end{prop}

Given a path $\gamma:[a, b]\rightarrow {\mathbb R}^n,$ we use the
notation
$$|\gamma|:=\{x\in {\mathbb R}^n:\,\exists\,t\in [a, b]: \gamma(t)=x\}$$
for the {\it locus} of $\gamma,$ see e.g. \cite[Section~1.1]{Va},
\cite[Section~II.1]{Ri}.

\medskip
{\it Proof of Theorem~\ref{th1}.}
In general, we follow the logic of the proof of Theorem~1.2 in~\cite{SevSkv}, see also Theorem~1.2
in~\cite{SSD} and Theorems~1--2 in~\cite{Sev$_3$}.
Let us fix $x, y\in K\subset {\mathbb B}^n$ and $f\in \mathfrak{F}_Q({\mathbb B}^n)$. We put
\begin{equation}\label{eq13A}
|f(x)-f(y)|:=\varepsilon_0\,.
\end{equation}
If $\varepsilon_0=0,$ there is nothing to prove. Let
$\varepsilon_0>0.$
Let us give a straight line through the points $f(x)$ and $f(y)$:
$r=r(t)=f(x)+(f(x)-f(y))t,$ $-\infty<t<\infty$.
%(see Figure~\ref{fig2}).
%
%\begin{figure}[h]
%\centerline{\includegraphics[scale=0.4]{Graphic2C.eps}} \caption{To
%the proof of Theorem~\ref{th1}}\label{fig2}
%\end{figure}
%
Let $\gamma^1:[1, c)\rightarrow {\mathbb B}^n,$ $1<c\leqslant \infty$
be a maximum $f$-lifting of the ray $r=r(t),$ $t\geqslant 1,$ with
the origin at the point $x,$ which exists due to
Proposition~\ref{pr3}. Let us to prove that, the case
$\gamma^1(t)\rightarrow x_1\in {\mathbb B}^n$ as $t\rightarrow c$ is
impossible. Indeed, in this case, by Proposition~\ref{pr3} we obtain
that $c=\infty$ and $f(x_1)=\lim\limits_{t\rightarrow +\infty}r(t).$
Due to the openness of $f,$ $f(x_1)\in f({\mathbb B^n}),$ but on the
other hand, $f(x_1)=\infty$ by the definition of $r=r(t).$ Since
$\infty\not\in f({\mathbb B^n}),$ we obtain a contradiction. Therefore,
$\gamma^1(t)\rightarrow x_1\in {\mathbb B}^n$ as $t\rightarrow c$, is impossible, as
required. By Proposition~\ref{pr3}
\begin{equation}\label{eq3B}
h(\gamma^1(t), \partial {\mathbb B}^n)\rightarrow 0
\end{equation}
as $t\rightarrow c-0.$ Similarly, denote by $\gamma^2:(d,
0]\rightarrow {\mathbb B}^n,$ $-\infty\leqslant d<0,$ the maximal
$f$-lifting of a ray $r=r(t),$ $t\leqslant 0,$ with the end at the
point $y,$ which exists by Proposition~\ref{pr3}. Similarly
to~(\ref{eq3B}), we obtain that
\begin{equation*}\label{eq4B}
h(\gamma^2(t), \partial {\mathbb B}^n)\rightarrow 0
\end{equation*}
as $t\rightarrow d+0.$ Let $z=\gamma^1(t_1)$ be some point on $\gamma^1,$ which is located at the distance $r_0/2$ from the unit
sphere, where $r_0:= d(K,\partial {\mathbb B}^n)$ and let $w=\gamma^2(t_2)$ be some point on
$\gamma^2,$ located at the distance $r_0/2$ from the unit sphere.
Put $\gamma^*:=\gamma^1|_{[1, t_1]}$ and $\gamma_*:=\gamma^1|_{[t_2,
0]}$. By the triangle inequality, ${\rm diam}\,(|\gamma^*|)\geqslant
r_0/2$ and ${\rm diam}\,(|\gamma_*|)\geqslant r_0/2.$ Let
$\Gamma:=\Gamma(|\gamma^*|, |\gamma_*|, {\mathbb B}^n)$. Now, by using \cite[lemma~4.3]{Vu$_2$}, we obtain that
\begin{equation}\label{eq7D}
M(\Gamma)\geqslant (1/2)\cdot M(\Gamma(|\gamma^*|, |\gamma_*|, {\mathbb
R}^n))\,,
\end{equation}
and, on the other hand, by~\cite[Lemma~7.38]{Vu$_1$}
\begin{equation}\label{eq7B}
M(\Gamma(|\gamma^*|, |\gamma_*|, {\mathbb R}^n))\geqslant
c_n\cdot\log\left(1+\frac1m\right)\,,
\end{equation}
where $c_n>0$ is some constant depends on $n$ only and
$$
m=\frac{{\rm dist}(|\gamma^*|, |\gamma_*|)}{\min\{{\rm diam\,}(|\gamma^*|),
{\rm diam\,}(|\gamma_*|)\}}\,.
$$
Note that, ${\rm diam}\,(|\gamma^i|)=\sup\limits_{\omega,w\in|\gamma^i|}|\omega-w|\geqslant r_0/2,$ $i=1,2$. Then, by~(\ref{eq7D}) and~(\ref{eq7B}) and taking into account that ${\rm
dist}\,(|\gamma^*|, |\gamma_*|)\leqslant |x-y|,$ we obtain
\begin{equation}
\label{eq7C}
M(\Gamma)\geqslant \widetilde{c_n}\cdot \log\left(1+\frac{r_0}{2{\rm
dist}(|\gamma^*|, |\gamma_*|)}\right)\geqslant \widetilde{c_n}\cdot
\log\left(1+\frac{r_0}{2|x-y|}\right)\,,
\end{equation}
where $\widetilde{c_n}>0$ is some constant depends on~$n$ only.
By the H\"{o}lder inequality, for any function $\rho\in {\rm adm}\,\Gamma$ we have
\begin{equation}\label{eq1B}
M(\Gamma)\leqslant \int\limits_{{\mathbb B}^n}\rho^n(x)\,dm(x)\leqslant
\left(\int\limits_{{\mathbb
B}^n}\rho^q(x)\,dm(x)\right)^{\frac{n}{q}}\cdot
(\Omega_n)^{\frac{q-n}{n}}\,.
\end{equation}
Taking in the right side of the inequality (\ref{eq1B}) the infimum over all $\rho\in {\rm
adm}\,\Gamma,$ we obtain that
\begin{equation}\label{eq1C} M(\Gamma)\leqslant {\rm inf}
\int\limits_{{\mathbb B}^n}\rho^n(x)\,dm(x)\leqslant
\left(M_q(\Gamma)\right)^{\frac{n}{q}}\cdot
(\Omega_n)^{\frac{q-n}{n}}\,.
\end{equation}
Combining~(\ref{eq7C}) and~(\ref{eq1B}), we obtain that
\begin{equation}\label{eq1D}
M_q(\Gamma)\geqslant
(m(\Omega_n))^{(n-q)q}{(\widetilde{c_n})}^{\frac{q}{n}}\cdot
\log^{\frac{q}{n}}\left(1+\frac{r_0}{2|x-y|}\right)\,.
\end{equation}

\medskip
Let $z_1:=f(z),$ $\varepsilon^{(1)}:=|f(x)- z^1|$ and
$\varepsilon^{(2)}:=|f(y)-z^1|.$ Note that
$$|f(y)- f(x)|+\varepsilon^{(1)}=$$
\begin{equation}\label{eq5B}
=|f(y)- f(x)|+|f(x)-z^1|= |z^1-f(y)|=\varepsilon^{(2)}\,,
\end{equation}
therefore, $\varepsilon^{(1)}<\varepsilon^{(2)}$.

Now let us to obtain an upper estimate for $M_q(\Gamma)$.
We put $\textbf{P}=|f(\gamma^*)|$, $\textbf{Q}=|f(\gamma^2)|$, and
$$A:=A(z^1, \varepsilon^{(1)}, \varepsilon^{(2)})=
\{x\in {\mathbb R}^n:
\varepsilon^{(1)}<|x-z^1|<\varepsilon^{(2)}\}\,.$$
Note that, $E:=\gamma^{\,*}$ and $F:=\gamma_{\,*}$ are continua in
${\mathbb B}^n.$ Let us to prove that
$$|\gamma^{\,*}|\subset
f^{\,-1}(\overline{B(z^1, \varepsilon^{(1)})})\,,\qquad
|\gamma_{\,*}|\subset f^{\,-1}(f({\mathbb B}^n)\setminus B(z^1,
\varepsilon^{(2)}))\,.$$
Indeed, let $x_*\in |\gamma^{\,*}|.$ Then $f(x_*)\in \textbf{P}$,
therefore, there exists numbers $1\leqslant t\leqslant s$ such that
$f(x_*)=f(y)+(f(x)-f(y))t,$ where $z^1=f(y)+(f(x)-f(y))s.$ Thus,
\begin{multline}
\label{eq2B}
|f(x_*)-z^1|=|(f(x)-f(y))(s-t)|
\\
\leqslant |(f(x)-f(y))(s-1)|=|(f(x)-f(y))s+f(y)-f(x))|
\\
=|f(x)-z^1|=\varepsilon^{(1)}\,.
\end{multline}
By~(\ref{eq2B}) it follows, that $|\gamma^{\,*}|\subset
f^{\,-1}(\overline{B(z^1, \varepsilon^{(1)})})$. The inclusion
$|\gamma_{\,*}|\subset f^{\,-1}(f({\mathbb B}^n)\setminus B(z^1,
\varepsilon^{(2)}))$ may be proved similarly.

\medskip
Let us put
$$\eta(t)= \left\{
\begin{array}{rr}
\frac{1}{\varepsilon_0}, & t\in [\varepsilon^{(1)}, \varepsilon^{(2)}],\\
0, & t\not\in [\varepsilon^{(1)}, \varepsilon^{(2)}]\,,
\end{array}
\right. $$
where $\varepsilon_0$ is a number from~(\ref{eq13A}). Note that the
function $\eta$ satisfies the relation~(\ref{eqA2}) for
$r_1=\varepsilon^{(1)}$ and $r_2=\varepsilon^{(2)}.$ Indeed,
by~(\ref{eq13A}) and (\ref{eq5B}) we obtain that
$$r_1-r_2=\varepsilon^{(2)}-\varepsilon^{(1)}=|f(y)-z^1|-|f(x)-
z^1|=$$$$=|f(x)-f(y)|=\varepsilon_0\,.$$
Then
$\int\limits_{\varepsilon^{(1)}}^{\varepsilon^{(2)}}\eta(t)\,dt=(1/\varepsilon_0)\cdot
(\varepsilon^{(2)}-\varepsilon^{(1)})\geqslant 1.$ Applying the
moduli inequality~(\ref{eq2*A}) for the point~$z^1,$ we obtain that
\begin{equation}\label{eq14***}
M_q(\Gamma)\leqslant \frac{1}{\varepsilon_0^q}\int\limits_{{\mathbb
R}^n} Q(z)\,dm(z)=\frac{\Vert Q\Vert_1}{{|f(x)-f(y)|}^{q}}\,.
\end{equation}
Finally, from (\ref{eq1D}) and (\ref{eq14***}) we obtain that
$$
(\Omega_n)^{(n-q)q}{(\widetilde{c_n})}^{\frac{q}{n}}\cdot
\log^{\frac{q}{n}}\left(1+\frac{r_0}{2|x-y|}\right)\leqslant
\frac{\Vert Q\Vert_1}{{|f(x)-f(y)|}^{q}}\,.
$$
Hence, It follows, that
$$|f(x)-f(y)|\leqslant C_n \cdot\frac{(\Vert Q\Vert_1)^{\frac{1}{q}}}
{\log^{\frac{1}{n}}\left(1+\frac{r_0}{2|x-y|}\right)}\,,$$
where
$C_n:=(\Omega_n)^{\frac{(q-n)q}{q}}{(\widetilde{c_n})}^{-\frac{1}{n}}$. The theorem is proved.~$\Box$

\section{H\"{o}lder continuity in arbitrary domains} 

Let $D, D^{\,\prime}$ be domains in ${\mathbb R}^n$, $n\geqslant 2$. For numbers $1\leqslant q<\infty$ and a Lebesgue measurable function $Q:{\mathbb
R}^n\rightarrow [0, \infty]$, $Q=0$ a.e. on $\mathbb R^n\setminus D^{\,\prime}$, we denote be $\mathfrak{R}^{q}_Q(D, D^{\,\prime})$ the family of all open and discrete mappings $f:D\rightarrow D^{\,\prime}$ such that the
moduli inequality~(\ref{eq2*A}) holds at any point $y_0\in D^{\,\prime}$. The following theorem generalizes \cite[Theorem~4.1]{SevSkv}.

\medskip
\begin{thm}\label{th6}
{\sl Let $Q\in L^1({\mathbb R}^n)$ and $q\geqslant n$. Suppose that, $K$ is compact in $D,$ and $D^{\,\prime}$ is bounded. Then there exists a constant $C=C(n,q, K, \Vert Q\Vert_1, D,D^{\,\prime})>0$ such that the inequality
\begin{equation}\label{eq2E}
|f(x)-f(y)|\leqslant C_n \cdot\frac{(\Vert Q\Vert_1)^{\frac{1}{q}}}
{\log^{\frac{1}{n}}\left(1+\frac{r_0}{2|x-y|}\right)}\,,\,\,r_0=d(K, \partial D),
\end{equation}
holds for any $x, y\in K$ and $f\in \mathfrak{R}_Q(D, D^{\,\prime})$, where $\Vert Q\Vert_1$ denotes the $L^1$-norm of the function $Q$ in ${\mathbb R}^n$. }
\end{thm}

\medskip
\begin{proof}
It is sufficient to find an upper bound for the value
\begin{equation}\label{eq1G}
|f(x)-f(y)|\cdot\log^{\frac{1}{n}}\left(1+\frac{r_0}{2|x-y|}\right)
\end{equation}
over all $x, y\in K$ and $f\in \mathfrak{R}_Q(D, D^{\,\prime}).$

\medskip
We fix $x, y\in K$ and $f\in \mathfrak{R}_Q(D, D^{\,\prime}).$ If
$|x-y|\geqslant r_0/2,$ the expression in~(\ref{eq1G}) is trivially
bounded. Indeed, by the triangle inequality,
\begin{equation}\label{eq16}
|f(x)-f(y)|\leqslant |f(x)|+|f(y)|\leqslant 2M_0\,,
\end{equation}
where $M_0=\sup\limits_{z\in D^{\,\prime}}|z|.$ Since $D^{\,\prime}$
is bounded, $M_0<\infty.$ By~(\ref{eq16}), we obtain that
\begin{equation}\label{eq2S}
|f(x)-f(y)|\cdot
\log^{\frac{1}{n}}\left(1+\frac{r_0}{2|x-y|}\right)\leqslant
M_0\cdot \log^{\frac{1}{n}}2\,,
\end{equation}
as required.

\medskip
Now let $|x-y|<r_0/2.$ In this case, $y\in B(x, r_0).$ Let $\psi$ be
a conformal mapping of the unit ball ${\mathbb B}^n$ onto the ball
$B(x, r_0),$ exactly, $\psi(z)=zr_0+x,$ $z\in {\mathbb B}^n.$ In
particular, $\psi^{\,-1}(B(x, r_0/2))=B(0, 1/2).$ Applying the
restriction $\widetilde{f}:=f|_{B(x, r_0)}$ and considering the
auxiliary mapping $F:=\widetilde{f}\circ \psi,$ $F:{\mathbb
B}^n\rightarrow D^{\,\prime},$ we conclude that the
relation~(\ref{eq2*A}) also holds for $F$ with the same
function~$Q.$ Then by Theorem~\ref{th1}
\begin{equation}\label{eq2F}
|F(\psi^{\,-1}(x))-F(\psi^{\,-1}(y))|\leqslant\frac{C_2\cdot (\Vert
Q\Vert_1)^{1/q}}{\log^{\frac{1}{n}}
\left(1+\frac{1}{4|\psi^{\,-1}(x)-\psi^{\,-1}(y)|}\right)}\,.
\end{equation}
Since $F(\psi^{\,-1}(x))=f(x)$ and $F(\psi^{\,-1}(y))=f(y),$ we may
rewrite~(\ref{eq2F}) in the form
\begin{equation}\label{eq2G}
|f(x)-f(y)|\leqslant\frac{C_2\cdot (\Vert
Q\Vert_1)^{1/q}}{\log^{\frac{1}{n}}
\left(1+\frac{1}{4|\psi^{\,-1}(x)-\psi^{\,-1}(y)|}\right)}\,.
\end{equation}
Note that, the mapping $\psi^{\,-1}(y)$ is Lipschitz with the
Lipschitz constant $\frac{1}{r_0}.$ In this case, due
to~(\ref{eq2G}), we obtain that
\begin{equation}\label{eq2H}
|f(x)-f(y)|\leqslant\frac{C_2\cdot (\Vert
Q\Vert_1)^{1/q}}{\log^{\frac{1}{n}}\left(1+\frac{r_0}{4|x-y|}\right)}\,.
\end{equation}
Finally, by the L'H\^{o}pital rule,
$\log^{\frac{1}{n}}\left(1+\frac{1}{nt}\right)\sim\log^{\frac{1}{n}}
\left(1+\frac{1}{kt}\right)$ as $t\rightarrow+0$ and any fixed $k,
n> 0.$ It follows that
$$\frac{C_2\cdot (\Vert
Q\Vert_1)^{1/q}}{\log^{\frac{1}{n}}\left(1+\frac{r_0}{4|x-y|}\right)}\leqslant
\frac{C_1\cdot (\Vert
Q\Vert_1)^{1/q}}{\log^{\frac{1}{n}}\left(1+\frac{r_0}{2|x-y|}\right)}$$
for some constant $C_1=C_1(r_0)> 0.$ Then, from~(\ref{eq2H}) it
follows that
\begin{equation}\label{eq2I}
|f(x)-f(y)|\leqslant\frac{C_1\cdot (\Vert
Q\Vert_1)^{1/q}}{\log^{\frac{1}{n}}\left(1+\frac{r_0}{2|x-y|}\right)}\,.
\end{equation}
Finally, from (\ref{eq2S}) and (\ref{eq2I}) it follows the desired
inequality~(\ref{eq2E}) with some constant
$$C:=\max\{C_1\cdot (\Vert Q\Vert_1)^{1/q}, M_0\cdot
\log^{\frac{1}{n}}2\}\,.~\Box$$
\end{proof}

\section{Boundary behavior of mappings} 

The following result  in the case $q=n$ was proved in \cite[Theorem~3.1]{SSD}, \cite[Theorem~4]{Sev$_3$}.

\medskip
\begin{thm}\label{th3}
{\sl\, Let $n\leqslant q<\infty,$ $D\subset {\mathbb R}^n,$
$n\geqslant 2,$ be a bounded domain with a weakly flat boundary, and
let $D^{\,\prime}\subset {\mathbb R}^n$ be a domain which is
finitely connected on its boundary. Suppose $f$ is open discrete and
closed mapping of $D$ onto $D^{\,\prime}$ satisfying the
relation~(\ref{eq2*A}) at any point $y_0\in
\partial D^{\,\prime},$ and the following condition holds:
for any $y_0\in \partial D^{\,\prime}$ and
$0<r_1<r_2<r_0:=\sup\limits_{y\in D^{\,\prime}}|y-y_0|$ there is
some set $E\subset[r_1, r_2]$ of positive linear Lebesgue measure
such that the function $Q$ is integrable on $S(y_0, r)$ for each
$r\in E.$ Then $f$ has a continuous extension
$\overline{f}:\overline{D}\rightarrow\overline{D^{\,\prime}},$
moreover, $\overline{f}(\overline{D})=\overline{D^{\,\prime}}.$

In particular, the statement of the theorem~\ref{th3} holds if $Q\in
L^1(D^{\,\prime}).$}
\end{thm}

\medskip
\begin{proof}
Let $x_0\in\partial D$.  We should prove the possibility of continuous extension of mapping $f$ to point $x_0$. Let us prove it from the opposite, namely, suppose that $f$
does not have a continuous extension to $x_0.$ Then, there are
sequences $x_i, y_i\in D,$ $i=1,2,\ldots ,$ such that $x_i,
y_i\rightarrow x_0$ as $i\rightarrow\infty,$ and, there is $a>0$
such that
\begin{equation}\label{eq9C}
h(f(x_i), f(y_i))\geqslant a>0
\end{equation}
for any $i\in {\mathbb N},$ where $h$ is a chordal (spherical) metric,
defined in~(\ref{eq3C}). Since the space $\overline{{\mathbb R}^n}$ is
compact, we may assume that $f(x_i)$ and $f(y_i)$ converge as
$i\rightarrow\infty$ to $z_1$ and $z_2,$ respectively, and
$z_1\ne\infty.$

Since $f$ is closed, it preserves the boundary of the domain
see~\cite[theorem~3.3]{Vu}, therefore $z_1, z_2\in
\partial D^{\,\prime}.$ Since $D^{\,\prime}$ is finitely connected
on its boundary, there are paths $\alpha:[0, 1)\rightarrow
D^{\,\prime}$ and $\beta:[0, 1)\rightarrow D^{\,\prime}$ such that
$\alpha\rightarrow z_1$ and $\beta\rightarrow z_2$ as $t\rightarrow
1-0$ such that $|\alpha|$ contains some subsequence of the sequence
$f(x_i)$ and $\beta$ contains some subsequence of the sequence
$f(y_i),$ $i=1,2,\ldots $ (see \cite[lemma~3.10]{Vu}). Without loss
of generalization, we may assume that the paths $\alpha$ and $\beta$
contain sequences $f(x_i)$ and $f(y_i),$ respectively.
%see Figure~\ref{fig1}.
%
%\begin{figure}[h]
%\centerline{\includegraphics[scale=0.5]{Graphic1C.eps}} \caption{To
%the proof of Theorem~\ref{th3}}\label{fig1}
%\end{figure}
%
%
Due to the definition of finite connectedness of the domain
$D^{\,\prime}$ on the boundary, we may assume that
\begin{equation}\label{eq8}
|\alpha|\subset B(z_1, R_*),\, |\beta|\subset {\mathbb R}^n\setminus
B(z_1, R_0)\,,\quad 0<R_*<R_0<\infty\,.
\end{equation}
We denote by $\alpha_i$ a subpath of $\alpha$ with the origin at a
point $f(x_i)$ and ends at $f(x_1)$ and, similarly, by $\beta_i$ a
subpath of $\beta$ starting at $f(y_i)$ and ending at $f(y_1).$ By
the change of a parameter, we may consider that, the paths
$\alpha_i$ and $\beta_i$ are parameterized so that $\alpha_i:[0,
1]\rightarrow D^{\,\prime}$ and $\beta_i:[0, 1]\rightarrow
D^{\,\prime}.$ Let $\widetilde{\alpha_i}:[0, 1)\rightarrow D$ and
$\widetilde{\beta_i}:[0, 1)\rightarrow D$ be whole $f$-liftings of
$\alpha_i$ and $\beta_i$ starting at points $x_i$ and $y_i,$
respectively (these lifts exist by~\cite[lemma~3.7]{Vu}). By
Proposition~\ref{pr3}, paths $\widetilde{\alpha_i}$ and
$\widetilde{\beta_i}$ can be extended to closed paths
$\widetilde{\alpha_i}:[0, 1]\rightarrow D$ and
$\widetilde{\beta_i}:[0, 1]\rightarrow D.$ Note that, the points
$f(x_1)$ and $f(y_1)$ may not have more than a finite number of
preimages under $f$ in $D,$ see~\cite[Theorem~2.8]{MS}. Then, there
is $r_0>0$ such that $\widetilde{\alpha_i}(1),
\widetilde{\beta_i}(1)\in D\setminus B(x_0, r_0)$ for all
$i=1,2,\ldots .$ Since the boundary of the domain $D$ is weakly
flat, for any $P>0$ there exists $i=i_P\geqslant 1$ such that
\begin{equation}\label{eq7}
M(\Gamma(|\widetilde{\alpha_i}|, |\widetilde{\beta_i}|,
D))>P\qquad\forall\,\,i\geqslant i_P\,.
\end{equation}
By H\"{o}lder inequality, for any function $\rho\in {\rm
adm}\,\Gamma,$
\begin{equation}\label{eq1BA}
M(\Gamma)\leqslant \int\limits_{D}\rho^n(x)\,dm(x)\leqslant
\left(\int\limits_{D}\rho^q(x)\,dm(x)\right)^{\frac{n}{q}}\cdot
m^{\frac{q-n}{n}}(D)\,.
\end{equation}
Letting~(\ref{eq1BA}) to $\inf$ over all $\rho\in {\rm
adm}\,\Gamma,$ we obtain that
\begin{equation}\label{eq1CC} M(\Gamma)\leqslant
\int\limits_{D}\rho^n(x)\,dm(x)\leqslant
\left(M_q(\Gamma)\right)^{\frac{n}{q}}\cdot m^{\frac{q-n}{n}}(D)\,.
\end{equation}
Using~(\ref{eq7}) and~(\ref{eq1CC}), we obtain that
\begin{equation}\label{eq7A}
M_q(\Gamma(|\widetilde{\alpha_i}|, |\widetilde{\beta_i}|, D))>P\cdot
m^{-\frac{q-n}{n}}(D)\qquad\forall\,\,i\geqslant i_P\,.
\end{equation}
Let us to show that, the condition~(\ref{eq7}) contradicts the
definition of mapping $f$ in~(\ref{eq2*A}). Indeed,
using~(\ref{eq8}) and applying~(\ref{eq2*A}) for
$E=|\widetilde{\alpha_i}|,$ $F=|\widetilde{\beta_i}|,$ $r_1=R_*$ and
$r_2=R_0,$ we obtain that
\begin{equation}\label{eq8B}
M_q(\Gamma(|\widetilde{\alpha_i}|, |\widetilde{\beta_i}|,
D))\leqslant \int\limits_{A(z_1,R_*,R_0)\cap D^{\,\prime}} Q(y)\cdot
\eta^{\,q}(|y-z_1|)\, dm(y)\,,
\end{equation}
where $\eta: (R_*,R_0)\rightarrow [0,\infty ]$ is any Lebesgue
measurable function such that
\begin{equation}\label{eq9B}
\int\limits_{R_*}^{R_0}\eta(r)\, dr\geqslant 1\,.
\end{equation}
Below, we use the standard conventions: $a/\infty=0$ for
$a\ne\infty,$ $a/0=\infty$ for $a>0$ and $0\cdot\infty=0$ (see,
e.g., \cite[3.I]{Sa}). Let us put $\widetilde{Q}(y)=\max\{Q(y),
1\},$
\begin{equation}\label{eq5V}
\widetilde{q}_{y_0}(r)=\frac{1}{\omega_{n-1}r^{n-1}}\int\limits_{S(y_0,
r)}\widetilde{Q}(y)\,d\mathcal{H}^{n-1}(y)
\end{equation}
and
\begin{equation}\label{eq13C}
I=\int\limits_{R_*}^{R_0}\frac{dt}{t^{\frac{n-1}{q-1}}
\widetilde{q}_{z_1}^{1/(q-1)}(t)}\,.
\end{equation}
By assumption of the theorem, for any $y_0\in \partial D^{\,\prime}$
and $0<r_1<r_2<r_0:=\sup\limits_{y\in D^{\,\prime}}|y-y_0|$ there is
a set $E\subset[r_1, r_2]$ of a positive Lebesgue linear measure
such that $Q$ is integrable on $S(y_0, r)$ for any $r\in E.$ Then
$0\ne I\ne \infty.$
In this case, the function
$\eta_0(t)=\frac{1}{It^{\frac{n-1}{q-1}}\widetilde{q}_{z_1}^{1/(q-1)}(t
)}$ satisfies the relation~(\ref{eq9B}). Substituting this function
into the right-hand side of~(\ref{eq8B}) and applying Fubini theorem
(see \cite[theorem~8.1, Ch.~III]{Sa}), we obtain that
$$M_q(\Gamma(|\widetilde{\alpha_i}|, |\widetilde{\beta_i}|,
D))\leqslant \int\limits_{A(z_1,R_*,R_0)\cap D^{\,\prime}} Q(y)\cdot
\eta^{\,q}(|y-z_1|)\, dm(y)=$$
\begin{equation}\label{eq14A}=
\int\limits_{R_*}^{R_0}\int\limits_{S(z_1, t)} Q(y)\cdot
\eta^{\,q}(|y-z_1|)\,d\mathcal{H}^{n-1}\,dt=
\frac{\omega_{n-1}}{I^{q-1}}<\infty\,.
\end{equation}
The relation~(\ref{eq14A}) contradicts~(\ref{eq7}), which disproves
the assumption made in~(\ref{eq9C}). The resulting contradiction
disproves the assumption that there is no a limit of $f$ at the
point $x_0.$

It remains to check the equality
$\overline{f}(\overline{D})=\overline{D^{\,\prime}}.$ It is obvious
that $\overline{f}(\overline{D})\subset\overline{D^{\,\prime}}.$ Let
us show that $\overline{D^{\,\prime}}\subset
\overline{f}(\overline{D}).$ Indeed, let $y_0\in
\overline{D^{\,\prime}},$ then either $y_0\in D^{\,\prime},$ or
$y_0\in \partial D^{\,\prime}.$ If $y_0\in D^{\,\prime},$ then
$y_0=f(x_0)$ and $y_0\in \overline{f}(\overline{D}),$ since by
condition $f$ is the mapping of $D$ onto $D^{\,\prime}.$ Finally,
let $y_0\in \partial D^{\,\prime},$ then there is a sequence $y_k\in
D^{\,\prime}$ such that $y_k=f(x_k)\rightarrow y_0$ as
$k\rightarrow\infty,$ $x_k\in D.$ Due to the compactness of
$\overline{{\mathbb R}^n},$ we may assume that $x_k\rightarrow x_0,$
where $x_0\in\overline{D}.$ Note that, $x_0\in \partial D,$ since
$f$ is open. Then $f(x_0)=y_0\in \overline{f}(\partial D)\subset
\overline{f}(\overline{D}).$ In the whole, Theorem~\ref{th3} is
proved, excluding the discussion of the situation $Q\in
L^1(D^{\,\prime}).$

\medskip
If $Q\in L^1(D^{\,\prime}),$ by the Fubini theorem,
$$\int\limits_{B(y_0, r_0)}Q(y)\, dm(y)=
\int\limits_{0}^{r_0}\int\limits_{S(y_0, t)}
Q(y)\,d\mathcal{H}^{n-1}\,dt<\infty\,,
$$
whence it follows that $q_{y_0}(t)<\infty$ for all $y_0\in\partial
D^{\,\prime}$ and almost all $t\in {\mathbb R}$ (here, of course, we
extend the function $Q$ by an identical zero outside
$D^{\,\prime}$). Thus, the case $Q\in L^1(D^{\,\prime})$ is a
special case of the conditions on $Q$ mentioned above. The theorem
is completely proved.~$\Box$
\end{proof}

\section{The equicontinuity of some family of mappings in the closure of domains}

{\it Proof of Theorem~\ref{th2}.} Let $f\in {\mathfrak
S}^{q}_{\delta, A, Q }(D, D^{\,\prime}).$ By Theorem~\ref{th3}, $f$
has a continuous extension $\overline{f}:\overline{D}\rightarrow
\overline{D^{\,\prime}},$ moreover,
$\overline{f}(\overline{D})=\overline{D^{\,\prime}}.$ The
equicontinuity of the family ${\mathfrak S}^{q}_{\delta, A, Q
}(\overline{D}, \overline{D^{\,\prime}})$ in $D$ is a statement of
Theorem~\ref{th6}. It remains to establish its equicontinuity on
$\partial D.$

We will carry out a proof from the opposite
(cf.~\cite[Theorem~1.2]{SSD}, \cite[Theorem~5]{Sev$_3$}). Assume
that, there is $x_0\in \partial D,$ a number $\varepsilon_0>0,$ a
sequence $x_m\in \overline{D},$ which converges to $x_0$ as
$m\rightarrow \infty$ and a sequence of mappings $\overline{f}_m\in
{\mathfrak S}^{q}_{\delta, A, Q }(\overline{D}, \overline{D})$ such that
\begin{equation}\label{eq12}
h(\overline{f}_m(x_m),\overline{f}_m(x_0))\geqslant\varepsilon_0,\quad
m=1,2,\ldots .
\end{equation}
Let us put $f_m:=\overline{f}_m|_{D}.$ Since $f_m$ has a continuous
extension on $\partial D,$ we may assume that $x_m\in D.$ Therefore,
$\overline{f}_m(x_m)=f_m(x_m).$ In addition, there exists a sequence
$x^{\,\prime}_m\in D$ such that $x^{\,\prime}_m\rightarrow x_0$ as
$m\rightarrow\infty$ and
$h(f_m(x^{\,\prime}_m),\overline{f}_m(x_0))\rightarrow 0$ as
$m\rightarrow\infty.$ Since the space $\overline{{\mathbb R}^n}$ is
compact, we may assume that the sequences $f_m(x_m)$ and
$\overline{f}_m(x_0)$ converge as $m\rightarrow\infty.$ Let
$f_m(x_m)\rightarrow \overline{x_1}$ and
$\overline{f}_m(x_0)\rightarrow \overline{x_2}$ as
$m\rightarrow\infty.$ By the continuity of the metric
in~(\ref{eq12}), $\overline{x_1}\ne \overline{x_2}.$ Since $f_m$ is
closed, it preserves the boundary (see~\cite[theorem~3.3]{Vu}). It
follows that $\overline{x_2}\in\partial D^{\,\prime}.$ Let
$\widetilde{x_1}$ and $\widetilde{x_2}$ be arbitrary distinct points
of the continuum $A,$ none of which coincides with $\overline{x_1}.$
Due to~\cite[Lemma~2.1]{SevSkv$_2$}, we may join two pairs of points
$\widetilde{x_1},$ $\overline{x_1}$ and $\widetilde{x_2},$
$\overline{x_2}$ using paths $\gamma_1:[0, 1]\rightarrow
\overline{D^{\,\prime}}$ and $\gamma_2:[0, 1]\rightarrow
\overline{D^{\,\prime}}$ such that $|\gamma_1|\cap
|\gamma_2|=\varnothing,$ $\gamma_1(t), \gamma_2(t)\in D$ for $t\in
(0, 1),$ $\gamma_1(0)=\widetilde{x_1},$
$\gamma_1(1)=\overline{x_1},$ $\gamma_2(0)=\widetilde{x_2}$ and
$\gamma_2(1)=\overline{x_2}.$ Since $D^{\,\prime}$ is locally
connected on $\partial D^{\,\prime},$ there are disjoint
neighborhoods $U_1$ and $U_2$ containing the points $\overline{x_1}$
and $\overline{x_2},$ such that the sets $W_i:=D^{\,\prime}\cap U_i$
are path connected. Without loss of generalization, we may assume
that $\overline{U_1}\subset B(\overline{x_1}, \delta_0)$ and
\begin{equation}\label{eq12C}
\overline{B(\overline{x_1},
\delta_0)}\cap|\gamma_2|=\varnothing=\overline{U_2}\cap|\gamma_1|\,,
\quad \overline{B(\overline{x_1}, \delta_0)}\cap
\overline{U_2}=\varnothing\,.
\end{equation}
Due to~(\ref{eq12C}), there is $\sigma_0>\delta_0>0$ such that
$$
\overline{B(\overline{x_1},
\sigma_0)}\cap|\gamma_2|=\varnothing=\overline{U_2}\cap|\gamma_1|\,,
\quad \overline{B(\overline{x_1}, \sigma_0)}\cap
\overline{U_2}=\varnothing\,.$$
We also may assume that $f_m(x_m)\in W_1$ and
$f_m(x^{\,\prime}_m)\in W_2$ for all $m\in {\mathbb N}.$ Let $a_1$ and
$a_2$ be two different points belonging to $|\gamma_1|\cap W_1$ and
$|\gamma_2|\cap W_2,$ in addition, let $0<t_1, t_2<1$ be such that
$\gamma_1(t_1)=a_1$ and $\gamma_2(t_2)=a_2.$ Join the points $a_1$
and $f_m(x_m)$ with a path $\alpha_m:[t_1, 1]\rightarrow W_1$ such
that $\alpha_m(t_1)=a_1$ and $\alpha_m(1)=f_m(x_m).$ Similarly, let
us join $a_2$ and $f_m(x^{\,\prime}_m)$ by a path $\beta_m:[t_2,
1]\rightarrow W_2,$ such that $\beta_m(t_2)=a_2$ and
$\beta_m(1)=f_m(x^{\,\prime}_m)$.
%(see Figure~\ref{fig4}).
%%
%\begin{figure}[h]
%\centerline{\includegraphics[scale=0.55]{Graphic4A.eps}} \caption{To
%the proof of Theorem~\ref{th2}}\label{fig4}
%\end{figure}
%%
%
Set
$$C^1_m(t)=\quad\left\{
\begin{array}{rr}
\gamma_1(t), & t\in [0, t_1],\\
\alpha_m(t), & t\in [t_1, 1]\end{array} \right.\,,\qquad
C^2_m(t)=\quad\left\{
\begin{array}{rr}
\gamma_2(t), & t\in [0, t_2],\\
\beta_m(t), & t\in [t_2, 1]\end{array} \right.\,.$$
Let $D^1_m$ and $D^2_m$ be total $f_m$-liftings of the paths
$|C^1_m|$ and $|C^2_m|$ starting at points $x_m$ and
$x^{\,\prime}_m,$ respectively (such lifts exist
by~\cite[Lemma~3.7]{Vu}). In particular, under the condition
$h(f_m^{\,-1}(A),
\partial D)\geqslant~\delta>0,$ which is part of the definition of the class ${\mathfrak
S}^{q}_{\delta, A, Q }(D, D^{\,\prime}),$ the ends of $b_m^1$ and
$b_m^2$ of paths $D^1_m$ and $D^2_m,$ respectively, distant from
$\partial D$ at a distance not less than $\delta.$

Denote by $|C^1_m|$ and $|C^2_m|$ the loci of the paths $C^1_m$ and
$C^2_m,$ respectively. Let us put
$$l_0=\min\{{\rm dist}\,(|\gamma_1|,
|\gamma_2|), {\rm dist}\,(|\gamma_1|, U_2\setminus\{\infty\})\}$$
and consider the coverage $A_0:=\bigcup\limits_{x\in |\gamma_1|}B(x,
l_0/4)$ of the path $|\gamma_1|$ using balls. Since $|\gamma_1|$ is
a compact set, we may choose a finite number of indices $1\leqslant
N_0<\infty$ and corresponding points $z_1,\ldots, z_{N_0}\in
|\gamma_1|$ such that $|\gamma_1|\subset
B_0:=\bigcup\limits_{i=1}^{N_0}B(z_i, l_0/4).$ In this case,
$$|C^1_m|\subset U_1\cup |\gamma_1|\subset
B(\overline{x_1}, \delta_0)\cup \bigcup\limits_{i=1}^{N_0}B(z_i,
l_0/4)\,.$$
Let us put
$$D_{mi}=f_m^{\,-1}\left(|C^1_m|\cap \overline{B(z_i,
l_0/4)}\right)\,,\quad 1\leqslant i\leqslant N_0\,,$$$$
D_{m0}=f_m^{\,-1}\left(|C^1_m|\cap \overline{B(\overline{x_1},
\delta_0)}\right)\,,\quad i=0\,.$$
Since $f_m$ is a closed mapping, the pre-image of an arbitrary
compact set in $D^{\,\prime}$ is a compact set in $D$ (see, e.g.,
\cite[Theorem~3.3 (4)]{Vu}). Then, the sets $D_{mi}$ are compact in
$D,$ and by the definition, $D_{mi}\subset
f_m^{\,-1}(\overline{B(z_i, l_0/4)})$ for $i>0$ and $D_{m0}\subset
f_m^{\,-1}(\overline{B(\overline{x_1}, \delta_0)}).$

Let $\Gamma^{\,*}_m$ be the family of all paths joining $|D^1_m|$
and $|D^2_m|$ in $D,$ and let $\Gamma_{mi}$ be a subfamily of paths
$\gamma:[0, 1]\rightarrow D$ in $\Gamma_m$ such that
$f(\gamma(0))\in \overline{B(z_i, l_0/4)}$ for $1\leqslant
i\leqslant N_0$ and $f(\gamma(0))\in \overline{B(\overline{x_1},
\delta_0)}$ for $i=0.$ In this case,
\begin{equation}\label{eq10C}
\Gamma^{\,*}_m=\bigcup\limits_{i=0}^{N_0}\Gamma_{mi}\,,
\end{equation}
where $\Gamma_{mi}$ is a family of all paths $\gamma:[0,
1]\rightarrow D$ such that $\gamma(0)\in D_{mi}$ and $\gamma(1)\in
|D^2_m|,$ $0\leqslant i\leqslant N_0.$
Due to the definition of $l_0$ and $\sigma_0,$
$$|D^2_m|\subset f^{\,-1}_m\left(D^{\,\prime}\setminus\left(\bigcup\limits_{i=1}^{N_0}
B(z_i, l_0/2)\cup B(\overline{x_1}, \sigma)\right)\right)\,.$$
Then, we may apply the definition of the class of mappings
in~(\ref{eq2*A}) to any family $\Gamma_{mi}$. Let us put
$\widetilde{Q}(y)=\max\{Q(y), 1\}$ and
$$\widetilde{q}_{z_i}(r)=\frac{1}{\omega_{n-1}r^{n-1}}\int\limits_{S(z_i,
r)}\widetilde{Q}(y)\,d\mathcal{H}^{n-1}\,.$$ Note that,
$\widetilde{q}_{z_i}(r)\ne \infty$ for $r\in E\subset
[l_0/4,l_0/2],$ $m_1(E)>0$ (this follows from the condition of the
theorem).
Let us put
$$I_i=I_i(z_i,l_0/4,l_0/2)=\int\limits_{l_0/4}^{l_0/2}\
\frac{dr}{r^{\frac{n-1}{q-1}}\widetilde{q}_{z_i}^{\frac{1}{q-1}}(r)}\,
,\quad 1\leqslant i\leqslant N_0\,,$$
$$I_0=I_0(\overline{x_1},\delta_0,\sigma_0)=\int\limits_{\delta_0}^{\sigma_0}\
\frac{dr}{r^{\frac{n-1}{q-1}}\widetilde{q}_{\overline{x_1}}^{\frac{1}{q-1}}(r
)}\,.$$
Note that, $I_i\ne 0,$ since $\widetilde{q}_{z_i}(r)\ne \infty$ for
$r\in E\subset [l_0/4,l_0/2],$ $m_1(E)>0.$ In addition,
$I_i\ne\infty,$ $i=0, 1,2, \ldots, N_0.$ In this case, we put
$$\eta_i(r)=\begin{cases}
\frac{1}{I_ir^{\frac{n-1}{q-1}}\widetilde{q}_{z_i}^{\frac{1}{q-1}}(r)}\,
,&
r\in [l_0/4, l_0/2]\,,\\
0,& r\not\in [l_0/4, l_0/2]\,,
\end{cases}$$
$$\eta_0(r)=\begin{cases}
\frac{1}{I_0r^{\frac{n-1}{q-1}}\widetilde{q}_{\overline{x_1}}^{\frac{1}{q-1}}(r
)}\,,&
r\in [\delta_0, \sigma_0]\,,\\
0,& r\not\in [\delta_0, \sigma_0]\,.
\end{cases}$$
Note that, the functions $\eta_i$ and $\eta_0$ satisfy~(\ref{eqA2}).
Substituting these functions into the definition~(\ref{eq2*A}), and
using the Fubini theorem with a ratio~(\ref{eq10C}), we obtain that
\begin{equation}\label{eq7E}
M_q(\Gamma^{\,*}_m)\leqslant
\sum\limits_{i=0}^{N_0}M_q(\Gamma_{im})\leqslant
\sum\limits_{i=1}^{N_0}\frac{\omega_{n-1}}{I_i^{q-1}}+
\frac{\omega_{n-1}}{I_0^{p-1}}:=C_0\,, \quad m=1,2,\ldots\,.
\end{equation}
Let us to show that, the relation~(\ref{eq7E}) contradicts the weak
flatness of the boundary of the domain $D^{\,\prime}.$ Indeed, by
construction
$$h(|D^1_m|)\geqslant h(x_m, b_m^1) \geqslant
(1/2)\cdot h(f^{\,-1}_m(A), \partial D)>\delta/2\,,$$
\begin{equation}\label{eq14}
h(|D^2_m|)\geqslant h(x^{\,\prime}_m, b_m^2) \geqslant (1/2)\cdot
h(f^{\,-1}_m(A), \partial D)>\delta/2
\end{equation}
for any $m\geqslant M_0$ and some $M_0\in {\mathbb N}.$
Put $U:=B_h(x_0, r_0)=\{y\in \overline{{\mathbb R}^n}: h(y,
x_0)<r_0\},$ where $0<r_0<\delta/4$ and the number $\delta$ refers
to ratio~(\ref{eq14}). Note that, $|D^1_m|\cap U\ne\varnothing\ne
|D^1_m|\cap (D\setminus U)$ for any $m\in{\mathbb N},$ because $h(|
D^1_m|)\geqslant \delta/2$ and $x_m\in |D^1_m|,$ $x_m\rightarrow
x_0$ at $m\rightarrow\infty.$ Similarly, $|D^2_m|\cap
U\ne\varnothing\ne |D^2_m|\cap (D\setminus U).$ Since $|D^1_m|$ and
$|D^2_m|$ are continua, by Proposition~\ref{pr2}
\begin{equation}\label{eq8A}
|D^1_m|\cap \partial U\ne\varnothing, \quad |D^2_m|\cap
\partial U\ne\varnothing\,.
\end{equation}
Let $C_0$ be the number from the relation~(\ref{eq7E}). Since
$\partial D$ is weakly flat, for the number $P:=C_0\cdot
m^{\frac{q-n}{n}}(D)>0$ there is a neighborhood $V\subset U$ of the
point $x_0$ such that
\begin{equation}\label{eq9A}
M(\Gamma(E, F, D))>C_0\cdot m^{\frac{q-n}{n}}(D)
\end{equation}
for any continua $E, F\subset D$ such that $E\cap
\partial U\ne\varnothing\ne E\cap \partial V$ and $F\cap \partial
U\ne\varnothing\ne F\cap \partial V.$ Let us show that,
\begin{equation}\label{eq10A}
|D^1_m|\cap \partial V\ne\varnothing, \quad |D^2_m|\cap
\partial V\ne\varnothing\end{equation}
for sufficiently large $m\in {\mathbb N}.$ Indeed, $x_m\in |D^1_m|$ and
$x^{\,\prime}_m\in |D^2_m|,$ where $x_m, x^{\,\prime}_m\rightarrow
x_0\in V$ as $m\rightarrow\infty.$ In this case, $|D^1_m|\cap
V\ne\varnothing\ne |D^2_m|\cap V$ for sufficiently large $m\in {\mathbb
N}.$ Note that $h(V)\leqslant h(U)\leqslant 2r_0<\delta/2.$
By~(\ref{eq14}), $h(|D^1_m|)>\delta/2.$ Therefore, $|D^1_m|\cap
(D\setminus V)\ne\varnothing$ and, therefore, $|D^1_m|\cap\partial
V\ne\varnothing$ (see Proposition~\ref{pr2}). Similarly,
$h(V)\leqslant h(U)\leqslant 2r_0<\delta/2.$ It follows
from~(\ref{eq14}) that, $h(|D^2_m|)>\delta/2.$ Therefore,
$|D^2_m|\cap (D\setminus V)\ne\varnothing.$ By
Proposition~\ref{pr2}, we obtain that $|D^2_m|\cap\partial
V\ne\varnothing.$ Thus, the ratio~(\ref{eq10A}) is established.
Combining relations~(\ref{eq8A}), (\ref{eq9A}) and (\ref{eq10A}), we
obtain that $M(\Gamma^{\,*}_m)=M(\Gamma(|D^1_m|, |D^2_m|,
D))>C_0\cdot m^{\frac{q-n}{n}}(D).$ Finally, by the H\"{o}lder
inequality, taking into account the last condition, we obtain that
\begin{equation}\label{eq11A}
M_q(\Gamma^{\,*}_m)\geqslant C_0\cdot m^{\frac{q-n}{n}}(D)\cdot
m^{-\frac{q-n}{n}}(D)=C_0.
\end{equation}
The latter relation contradicts with~(\ref{eq7E}), which proves
theorem in the case of functions $Q$ integrable over spheres. The
case $Q\in L^1(D^{\,\prime})$ can be considered by analogy with the
last one part of the proof of Theorem~\ref{th3}.~$\Box$

\medskip
{\bf 6. Consequences for mappings with other modulus and capacity
conditions.} First of all, consider the relation
\begin{equation}\label{eq12A}
M_q(\Gamma(E, F, D))\leqslant \int\limits_{f(D)} Q(y)\cdot
\rho_*^{\,q}(y)\,dm(y)\qquad\forall\,\,\rho_*\in {\rm
adm}(f(\Gamma(E, F, D)))\,.
\end{equation}
The following statement holds.

\medskip
\begin{thm}\label{th7}
{\sl\, Let $y_0\in f(D),$ $q<\infty$ and let $Q:{\mathbb
R}^n\rightarrow [0, \infty]$ be a Lebesgue measurable function. If
$f$ is a mapping that satisfies relation~(\ref{eq12A}) for any
disjoint nondegenerate compact sets $E, F\subset D,$ then $f$ also
satisfies condition~(\ref{eq2*A}) for arbitrary compact sets
$E\subset f^{\,-1}(\overline{B(y_0, r_1)}),$ $F\subset
f^{\,-1}(f(D)\setminus B(y_0, r_2)),$
$0<r_1<r_2<r_0=\sup\limits_{y\in D^{\,\prime}}|y-y_0|,$ and an
arbitrary Lebesgue measurable function $\eta: (r_1,r_2)\rightarrow
[0,\infty ]$ with the condition~(\ref{eqA2}).}
\end{thm}

\medskip
\begin{proof}
Let $E\subset f^{\,-1}(\overline{B(y_0, r_1)}),$ $F\subset
f^{\,-1}(f(D)\setminus B(y_0, r_2)),$
$0<r_1<r_2<r_0=\sup\limits_{y\in D^{\,\prime}}|y-y_0|,$ be arbitrary
non-degenerate compacta. Also, let $\eta: (r_1,r_2)\rightarrow
[0,\infty]$ be an arbitrary Lebesgue measurable function that
satisfies condition~(\ref{eqA2}). Let us put
$\rho_*(y):=\eta(|y-y_0|)$ for $y\in A\cap f(D)$ and $\rho_*(y)=0$
otherwise, where $A=A(y_0, r_1, r_2)=\{y\,\in\,{\mathbb R}^n :
r_1<|y-y_0|<r_2\}.$ By Luzin theorem, we may assume that the
function $\rho_*$ is Borel measurable (see e.g.,
\cite[Section~2.3.6]{Fe}). By~\cite[theorem~5.7]{Va}
$$\int\limits_{\gamma_*}\rho_*(y)\,|dy|\geqslant
\int\limits_{r_1}^{r_2}\eta(r)\,dr\geqslant 1$$
for any (rectifiable) path $\gamma_*\in \Gamma(f(E), f(F), f(D)).$
Then, by~(\ref{eq12A}), we obtain that
$$M_q(\Gamma(E, F, D))\leqslant \int\limits_{A\cap f(D)}
Q(y)\cdot\rho_*^q(y)\,dm(y)=\int\limits_{A\cap f(D)} Q(y)\cdot
\eta^q(|y-y_0|)\, dm(y)\,.~\Box$$
\end{proof}

\medskip
Given a Lebesgue measurable function
$Q:{\mathbb R}^n\rightarrow [0, \infty],$ a {\it $q$-capacity of
$(E, F)$ with a weight $Q$ and with a respect to $D$} is defined by
\begin{equation}\label{eq5G}{\rm cap}_{q, Q}\,(E, F, D)\quad=\quad\inf\limits_{u\,\in\,W_0(E , F)}
\quad\int\limits_D\,Q(x)\cdot|\nabla u|^q\,\,dm(x)\,.
\end{equation}
The following statement holds.

\medskip
\begin{thm}\label{th8}
{\sl\, Let $y_0\in f(D),$ $q<\infty$ and let $Q:{\mathbb
R}^n\rightarrow [0, \infty]$ be Lebesgue measurable function. If $f$
is a homeomorphism that satisfies the relation
\begin{equation}\label{eq13B}
{\rm cap}_q(E, F, D)\leqslant {\rm cap}_{q, Q}\,(f(E), f(F),
f(D))\,,
\end{equation}
for arbitrary compacts (continua) $E, F\subset D,$ and
\begin{equation}\label{eq14B}{\rm cap}_{q, Q}\,(f(E), f(F), f(D))=M_{q, Q}\,(f (E), f(F),
f(D))\,,
\end{equation}
where
$$M_{q, Q}\,(f(E),
f(F), f(D))=\inf\limits_{\rho_*\in{\rm adm\,}\Gamma(f(E), f(F),
f(D))}\int\limits_{f(D)}\rho_*^q(y)\cdot Q(y)\,dm(y)\,,$$
then $f$ satisfies the condition~(\ref{eq2*A}) for arbitrary
compacts (continua) sets $E\subset f^{\,-1}(\overline{B(y_0,
r_1)}),$ $F\subset f^{\,-1}(f(D)\setminus B(y_0, r_2)),$
$0<r_1<r_2<r_0=\sup\limits_{y\in D^{\,\prime}}|y-y_0|,$ and an
arbitrary Lebesgue measurable function $\eta: (r_1,r_2)\rightarrow
[0,\infty ]$ with the condition~(\ref{eqA2}).}
\end{thm}

\medskip
\begin{proof}
Let $E\subset f^{\,-1}(\overline{B(y_0, r_1)}),$ $F\subset
f^{\,-1}(f(D)\setminus B(y_0, r_2)),$
$0<r_1<r_2<r_0=\sup\limits_{y\in D^{\,\prime}}|y-y_0|,$ be arbitrary
nondegenerate compacta. Also, let $\eta: (r_1,r_2)\rightarrow
[0,\infty]$ be an arbitrary Lebesgue measurable function that
satisfies the condition~(\ref{eqA2}). By Hesse equality (see
\cite[Theorem~5.5]{Hes}), ${\rm cap}_q(E, F, D)=M_q(\Gamma(E, F,
D)).$ Since $f$ is a homeomorphism, $f(\Gamma(F, E, D))=\Gamma(f(E),
f(F), f(D)).$ Then, by~(\ref{eq13B}) we obtain that
\begin{equation}\label{eq13D}
M_q(\Gamma(E, F, D))\leqslant{\rm cap}_{q, Q}\,(f(E), f(F),
f(D))\leqslant \int\limits_{f(D)} Q(y)\cdot \rho_*^{\,q}(y)\,dm(y)
\end{equation}
for any function $\rho_*\in {\rm adm}\,f(\Gamma(E, F, D))={\rm
adm}\,\Gamma(f(E), f(F), f(D)).$ The desired conclusion follows by
Theorem~\ref{th7}.~$\Box$
\end{proof}

\medskip
Due to Theorem~\ref{th8}, all results of this paper hold for
homeomorphisms with~(\ref{eq13D}), the corresponding weight $Q$ of
which satisfies the relation~(\ref{eq14B}).

\vskip 0.5cm

Vladimir Gol'dshtein; Department of Mathematics, Ben-Gurion University of the Negev, P.O.Box 653, Beer Sheva, 8410501, Israel

\emph{E-mail address:} \email{vladimir@math.bgu.ac.il} \\

Evgeny Sevost'yanov; Department of Mathematical Analysis, Zhytomyr
Ivan Franko State University, 40 Velyka Berdychivs'ka Str.,
Zhytomyr, 10008, Ukraine

Institute of Applied Mathematics and Mechanics of NAS of Ukraine, 19
Henerala Batyuka Str., 84 100 Slov'yans'k, Ukraine

\emph{E-mail address:} \email{esevostyanov2009@gmail.com} \\

Alexander Ukhlov; Department of Mathematics, Ben-Gurion University of the Negev, P.O.Box 653, Beer Sheva, 8410501, Israel

\emph{E-mail address:} \email{ukhlov@math.bgu.ac.il

\end{document}